\documentclass[10pt]{amsart}

\newcommand{\Rp}{(0,\IF)}
\newcommand{\dln}{\, d \upsilon(\x)}
\newcommand{\ir}{(0,\IF)^2}
 \def\x{x_1,x_2}

 \newcommand{\HR}{\notag \\}

\usepackage[square,compress,comma, numbers,sort]{natbib}
\usepackage[colorlinks=true, citecolor=red, linkcolor=blue]{hyperref}
\usepackage{amsfonts,mathtools}

\allowdisplaybreaks[4]

\usepackage{amssymb}
\usepackage{color}
\newcommand{\abs}[1]{\left\lvert #1 \right\rvert}
\newcommand{\Abs}[1]{ \biggl \lvert #1 \biggr \rvert}

\DeclarePairedDelimiterXPP\pk[1]{\mathbb{P}}\{ \}{}{ #1}
\DeclarePairedDelimiterXPP\E[1]{\mathbb{E}}\{ \}{}{	#1}

\DeclarePairedDelimiterXPP\ind[1]{\mathbb{I}}( ){}{	#1}

\def\FRE{\mbox{Fr\'{e}chet }}

 \usepackage{xparse}

\NewDocumentCommand{\ceil}{s O{} m}{%
  \IfBooleanTF{#1} 
    {\left\lceil#3\right\rceil} 
    {#2\lceil#3#2\rceil} 
}
\NewDocumentCommand{\floor}{s O{} m}{%
  \IfBooleanTF{#1} 
    {\left\lfloor#3\right\rfloor}
    {#2\lfloor#3#2\rfloor}
}

\definecolor{c20}{rgb}{0.,0.7,0.}
\definecolor{c30}{rgb}{0.,0.,1.}
\definecolor{c40}{rgb}{1,0.1,0.7}
\definecolor{c50}{rgb}{1,0,0}
\definecolor{c60}{rgb}{1,0.9,0.1}
\definecolor{c70}{rgb}{0.50,1.00,0.00}

\numberwithin{equation}{section}
\newtheorem{theo}{Theorem}[section]
\newtheorem{sat}[theo]{Proposition}
\newtheorem{de}[theo]{Definition}
\newtheorem{lem}[theo]{Lemma}

\newtheorem{example}[theo]{Example}
\newtheorem{korr}[theo]{Corollary}
\newtheorem{remark}[theo]{Remark}

\numberwithin{equation}{section}

\newcommand{\prooftheo}[1]{ \textsc{Proof of Theorem} \ref{#1} }
\newcommand{\proofprop}[1]{\textsc{Proof of Proposition} \ref{#1}}
\newcommand{\prooflem}[1]{\textsc{Proof of Lemma} \ref{#1}}
\newcommand{\proofkorr}[1]{\textsc{Proof of Corollary} \ref{#1}}

\newcommand{\QED}{\hfill $\Box$}

\newcommand{\COM}[1]{}

\def\IF{\infty}

\newcommand{\R}{\mathbb{R}}
\newcommand{\inr}{\in \R}

\newcommand{\BQN}{\begin{eqnarray}}
\newcommand{\EQN}{\end{eqnarray}}
\newcommand{\BQNY}{\begin{eqnarray*}}
\newcommand{\EQNY}{\end{eqnarray*}}
\def\ldot{, \ldots,}

\newcommand{\limit}[1]{\lim_{#1 \to   \infty}}

\newcommand{\kb}[1]{\boldsymbol{#1}}
\newcommand{\vk}[1]{\kb{#1}}

\def\bqny#1{\begin{eqnarray*} #1 \end{eqnarray*}}
\def\bqn#1{\begin{eqnarray} #1 \end{eqnarray}}

\newcommand{\BS}{\begin{sat}}
\newcommand{\ES}{\end{sat}}
\newcommand{\BT}{\begin{theo}}
\newcommand{\ET}{\end{theo}}
\newcommand{\BK}{\begin{korr}}
\newcommand{\EK}{\end{korr}}

\newcommand{\BEX}{\begin{example}}
\newcommand{\EEX}{\end{example}}

\newcommand{\BD}{\begin{de}}
\newcommand{\ED}{\end{de}}

\newcommand{\BIT}{\begin{itemize}}
\newcommand{\EIT}{\end{itemize}}
\newcommand{\BDI}{\begin{description}}
\newcommand{\EDI}{\end{description}}

\newcommand{\BRM}{\begin{remark}}
\newcommand{\ERM}{\end{remark}}

\newcommand{\BEL}{\begin{lem}}
\newcommand{\EEL}{\end{lem}}

\newcommand{\nelem}[1]{{Lemma \ref{#1}}}
\newcommand{\neprop}[1]{{Proposition \ref{#1}}}
\newcommand{\netheo}[1]{{Theorem \ref{#1}}}

\newcommand{\nekorr}[1]{{Corollary \ref{#1}}}

\def\TT{\mathcal{T}}
\def\TT{\R }

\def\Z{\mathbb{Z}}
\def\inn{\in \mathbb{N}}

\begin{document}

\title{Domination of Sample Maxima and Related Extremal  Dependence Measures}

\author{Enkelejd  Hashorva}
\address{Enkelejd Hashorva, Department of Actuarial Science 
University of Lausanne,\\
UNIL-Dorigny, 1015 Lausanne, Switzerland}
\email{Enkelejd.Hashorva@unil.ch}

\bigskip

 \maketitle

\begin{quote}

\def\ee{\mu(H,Q)}
{\bf Abstract}: For a given $d$-dimensional distribution function (df) $H$ we introduce the   class of dependence measures  $ \ee = - \E*{ \ln H(Z_1 \ldot Z_d)}, $ where the random vector $(Z_1 \ldot Z_d)$ has df   $Q$ which has the same marginal  df's   as $H$. If both $H$ and $Q$ are max-stable df's, we show that for a df  $F$ in the max-domain of attraction of $H$, this dependence measure explains the extremal dependence exhibited by $F$. Moreover we prove that $\ee$ is the limit of the probability that the  maxima of a random sample from $F$ is marginally dominated by some random vector with df in the max-domain of attraction of $Q$. We show a similar result for the complete domination of the sample maxima which leads to another measure of dependence denoted by $\lambda(Q,H)$. In the literature $\lambda(H,H)$ with $H$ a max-stable df  has been studied in the context of records,  multiple maxima, concomitants of order statistics and  concurrence probabilities. It turns out that both  $\ee$ and $\lambda(Q,H)$ are closely related. If $H$ is max-stable we derive useful representations for both $\ee$ and $\lambda(Q,H)$. Our applications include equivalent conditions for $H$ to be a product df   and $F$ to have asymptotically independent components. 
\end{quote}



\def\TTT{\mathcal{T}}
 \def\TT{\mathcal{T}}
  \def\intT{\int_{\TT}}
  \def\intTT{\sum_{t \in \delta \Z} }
  \def\intTTT{\sum_{t \in \Z} }

 \def\TT{T}

\section{Introduction} 
Let $H$ be a $d$-dimensional distribution function (df) with 
unit \FRE marginal df's  $\Phi(x)=e^{-1/x},x>0$.  We shall assume in the sequel that $H$ is a max-stable df,  which in our setup is equivalent with the homogeneity property 
\bqn{\label{hany} 
H^t(x_1 \ldot x_d)= H(t x_1 \ldot t x_d)
}
 for any $t>0, x_i \in (0,\IF), 1 \le i\le d $, see e.g., \cite{Res1987, Segers,Faletal2010}.
 The class of max-stable df's is very large with two extreme instances 
$$H_0(x_1 \ldot x_d)= \prod_{i=1}^d \Phi(x_i), \quad  H_\IF(x_1 \ldot x_d)= \min_{1 \le i \le d} \Phi(x_i)$$
  the product df  $H_0$ and the upper df $H_\IF$, respectively. Hereafter  $\bar G=1- G$ stands for the survival function of some univariate df $G$. It follows easily by the lower \FRE-Hoeﬀding  bound that 
  \bqn{ (H(nx_1 \ldot nx_d ))^n &\ge&  \Bigl(\max\big (0, 1- \sum_{i=1}^d \bar \Phi(n x_i))\Bigr)^n \ge  e^{ \liminf_{n\to \IF} n \ln  (1- \sum_{i=1}^d \bar \Phi(n x_i))}\HR 
  	&= & H_0(x_1 \ldot x_d), \quad x_i \in (0,\IF), \quad  i\le d. \label{ghr}
 }	 
Indeed, \eqref{ghr}  is well-known and follows for instance using the Pickands representation of $H$,  see e.g., \cite{Faletal2010}[Eq. (4.3.1)] or  the inf-argmax formula as shown  in Section 4. 
  Consequently, any max-stable df  $H$ lies between $H_0$ and $H_\IF$, i.e.,
\bqn{ \label{assoc} 
	H_0(x_1 \ldot x_d)\le  H(x_1\ldot x_d)\le H_\IF(x_1 \ldot x_d), \quad x_i \in (0,\IF), 1 \le i\le d.
}	
From multivariate extreme value theory, see e.g., \cite{deHaan, Res1987, Segers, Faletal2010} we know that  $d$-dimensional max-stable df's $H$ are limiting df's of the component-wise maxima of $d$-dimensional independent and identically distributed (iid) random vectors with some df  $F$. In that case, $F$ is said to be in the max-domain of attraction (MDA) of $H$ (abbreviated $F \in MDA(H)$). For simplicity we shall assume throughout in the following that $F$ is a df on $[0,\IF)^d$  with marginal df's $F_i\in MDA(\Phi),i\le d$ that have norming constants $a_n=n, n\inn$, and thus we have 
\bqn{ \label{nco} 
	\limit{n} F_i^n(nx)= \Phi(x), \quad x\in \R 
}	 for all $i\le d$, where we set $\Phi(x)=0$ if  $x\le 0$. Consequently,  $F$ is in the MDA of some max-stable df $H$ if  further
\bqn{ \label{plakut} \limit{n} \sup_{x_i \in \R, 1\le i\le d} \Abs{F^n(nx_1 \ldot nx_d)- H(x_1 \ldot x_d)}=0.
}
In the special case that $F$ has asymptotically  independent marginal df's, meaning that for $(X_1 \ldot X_d)$ with df $F$  
 \bqn{ \limit{n} n\pk{ X_i >n x_i , X_j > n x_j}=0, \quad x_i,x_j \in (0,\IF), \quad \forall i \not=j \le d,
\label{turbo}
}
then  $F\in MDA(H_0)$ if simply $F_i \in MDA( \Phi), i\le d$.\\
 In various applications it is important to be able to determine if some max-stable df $H$ resulting from the approximation in \eqref{plakut}  is equal to $H_0$, which in the light of multivariate extreme value theory means that the component-wise maxima $\vk{M}_n:=(\max_{1 \le i \le n} X_{i1}  \ldot \max_{1 \le i \le n} X_{id} ),n\ge 1$ of a $d$-dimensional random sample $(X_{i1} \ldot X_{id}), i=1 \ldot n$ of size $n$ from $F$ has asymptotically  independent  components. \\
The strength of dependence of the components of $\vk M_n$, or in other words the extremal dependence manifested in $F$, in view of the approximation \eqref{plakut} can be measured by calculating some appropriate dependence measure for $H$ (when the limiting df $H$ is known).\\
For any  random vector $\vk Z=(Z_1 \ldot Z_d)$ with df $Q$ which has the same marginal df's as $H$   we introduce a class of dependence measure for $H$ indexed by $Q$ given by
\def\ee{\mu(H,Q)}
\bqn{
	 \ee = - \E*{ \ln H(Z_1 \ldot Z_d) }. 
	}
In view of \eqref{assoc}, since $- \ln H_i(Z_i)$ is a unit exponential random variable, we have 
\bqn{\label{gabi}  
	1 =  \max_{1\le i \le d} \E{ - \ln H_i(Z_i)} \le  - \E[\Big]{\ln \min_{1 \le i\le d} H_i(Z_i)} \le  \ee \le - \E[\Big]{\ln \prod_{i=1}^d H_i(Z_i)} =  d
}
and in particular 
\bqn{\label{eqbd}
	 \mu(H_0,Q)=d,  \quad \mu(H_\IF,H_\IF)=1.
}
Clearly, $\ee$ can be defined for any df $H$ and it does not depend on the choice of the marginal df's of $H$. In this contribution we shall show that $\ee$ is particularly interesting for $H$ being max-stable.

\def\upi{\underline{\pi}_n}
\def\lpi{\overline{\pi}_n}

Next, consider the case that both $H$ and $Q$ are max-stable. It follows that (see  \netheo{thm20}) for $F$ satisfying \eqref{plakut} and $G \in MDA(Q)$   
\bqn{\label{eqJamba} \mu(Q,H) = \limit{n} \mu_n(G,F^n ), \quad \mu_n(G,F^n)=n\int_{\R^d} [1- G(x_1 \ldot x_d)] \, d F^n(x_1 \ldot x_d),
}
provided that both $F$ and $G$ are continuous. 
In view of \eqref{eqJamba}, we see that  $\ee$  relates to $F$ under \eqref{plakut}.\\
Let in the following $\vk W$ denote a random vector   with df $G$ being independent of $\vk M_n$. 
 We say that  $\vk W$  marginally dominates $\vk M_n$, if  
there exists some $i\le d$ such that $W_i > M_{ni}$. Consequently, 
assuming further  that $\vk W$ is independent of $\vk M_n$ we have  
$$ \frac{\mu_n(G,F^n)}{n}= \pk{ \vk W \text{ marginally dominates } \vk M_n}=: \upi.$$
Re-writing \eqref{eqJamba} we have  $ \limit{n} n\upi= \ee$ and thus  $\ee$ appears naturally in the context of marginal dominance of sample maxima.\\
Our motivation for introducing $\ee$ comes from results and ideas of A. Gnedin, see \cite{MR1240421,Gnedin94,Gnedin98} where multiple maxima of random samples is investigated. In the turn, the probability of observing a multiple maximum is closely related to the complete domination of sample maxima as we shall explain below.\\ We say that $\vk W$ completely dominates $\vk{M}_n$ if $W_i >   M_{ni}$ for any $i\le d$. Assuming that $F$ and $G$ are continuous, we have  
$$ 
\lambda_n(F^n, G):= n\int_{\R^d} F^n(x_1 \ldot x_d) \, d G(x_1 \ldot x_d)=
 n \pk{ \vk W \text{ completely  dominates } \vk M_n}=:n\lpi .$$
If further $F\in MDA(H), G\in MDA(Q)$  we show in \netheo{thm20} below  that
\bqn{\label{lamb}
	\limit{n} \lambda_n(G^n, F)= \lambda(Q,H)= \int_{(0,\IF)^d} Q(x_1 \ldot x_d) \, d \upsilon(x_1 \ldot x_d),
}
where $\upsilon$ denotes the exponent measure of $H$ defined on $E= [0,\IF]^d \setminus (0 \ldot 0)$, see \cite{Res1987,Faletal2010} for more details on the exponent measure. Note in passing that by symmetry $\limit{n}\lambda_n(F^n, G)= \lambda(H,Q)$ follows.

Our notation and definitions of $\lpi$ and $\upi$ agree with those in \cite{DombExtrC} for the particular case that $F=G$. Therein the complete and simple records are discussed.  If $F$ is continuous and  $F=G$ we have that $(n+1)\lpi$ equals 
$$\pk*{ \max_{1 \le i \le n+1} X_{ij} = X_{1j}, j=1 \ldot d },$$ which is the probability of observing a multiple maximum,    see  \cite{Gnedin94, gnedin1994conical, Gnedin98, hashorva2002b,hashorva2005,StoevDOMB}. 
  There are only few contributions that discuss the asymptotics of $\lambda_n(G^n, F)$ for $F\not=G$, see \cite{hashorva2001,MR2002123, hashorva2002}.

Since the exponent measure can be defined also for max-id.\ df $H$, i.e., if $H^t$ is a df for any $t>0$, then as above $\lambda(Q,H)$ can also be defined for any such df $H$ and any given $d$-dimensional df $Q$.  We shall show that $\ee$ and $\lambda(Q,H)$ are closely related. In particular, for $d=2$ we have $\ee= 2- \lambda(Q,H)$, provided that $H$ is a max-id.\ df. In particular, we show how to define $\lambda(Q,H)$ for any $H$ and $Q$.\\
For $H$ being a max-id.\ df we also show how to calculate $\ee$ and $\lambda(Q,H)$ by a limiting procedure, which relates to domination of $d$-dimensional random vectors, see \netheo{GER} below. 

It turns out that both dependence measures $\ee$ and $\lambda(Q,H)$ are very tractable if $H$ is max-stable (note that such $H$ is also max-id.\ df). In particular, we show that 
$\ee$ is the extremal coefficient of some $d$-dimensional max-stable df $H^*$, i.e., $\ee= - \ln H^*(1 \ldot 1)$.  Moreover, we derive in \netheo{kleefeld} tractable expressions for $\ee$ and $\lambda(Q,H)$, which are useful for simulations of these dependence measures if  the de Haan spectral representation of $H$ is known. 

It is of particular interest  for multivariate extreme value theory to derive tractable criteria that  identify if a max-stable df $H$ is equal to $H_0$. In our first application  we show several equivalent conditions to  $H=H_0$.  
  
In view of \eqref{eqJamba} and \eqref{lamb} we see that both measures of extremal dependence $\ee$ and $\lambda(Q,H)$  capture the extremal properties of $F \in MDA(H)$. Motivated by the relation between $\ee$ and $\lambda(Q,H)$ we derive in our second application several conditions equivalent to \eqref{turbo}. 

Both $\ee$ and $\lambda(Q,H)$ can be  defined for any $d$-dimensional df $H$ and $Q$. When $H$ is max-stable, these are dependence measures for $H$, since independent of the choice of $Q$, we can determine if $H=H_0$, see \neprop{thm1}, statement ii). A simple choice for $Q$ is taking $Q=H$. Alternatively, one can take $Q=H_0$ or $Q=H_\IF$. Independent of the choice of $Q$ we show in \neprop{thm1} that $\mu(H,Q)=2$ is equivalent with $H=H_0$. In particular, this result shows that $\mu(H,Q)$ is a measure of dependence of $H$ (and not for $Q$). 
 
 Brief organisation of the rest of the paper: In  Section 2 we derive the basic properties of both measures of $\ee$ and $\lambda(Q,H)$ if $H$ is a max-id.\ df. More tractable formulas are then derived for $H$ being a max-stable df. Section 3 is dedicated to applications. We present some auxiliary results in Section 4 followed by the proofs of the main results in Section 5.

\section{Main Results}  
In the following $H$ and $Q$ are  $d$-dimensional df's with unit \FRE marginals  df's and $\vk Z$ is a random vector with df   $Q$.  
 The second dependence measure $\lambda(Q,H)$ defined  in \eqref{lamb} is determined in terms of the exponent measure $\nu$ of $H$, under the max-stability assumption on $H$.\\
  A larger class of multivariate df's  is that of max-id.\ df's. Recall that $H$ is max-id., if $H^t$ is a df for an $t>0$. For such df's the corresponding exponent measure can be constructed, see e.g., \cite{Res1987}, and therefore we can define $\lambda(Q,H)$ as in the Introduction for any $H$ a max-id.\ df and any given df $Q$. Note that any max-stable df is a max-id. df, therefore in the following we shall consider first the general case that $H$ is a max-id.\ df, and then focus on the more tractable case that $H$ is a max-stable df.
  
\subsection{Max-id.\ df $H$}
Our analysis shows that    $\ee$ and $\lambda(Q,H)$ are closely related. Specifically,  
if $d=2$, then $\ee = 2- \lambda(Q,H)$, provided that $H$ is a max-id.\ df. Such a relationship does not hold for the case $d>2$. However as we show below it is possible to calculate $\ee$ if we know $\lambda(Q_K,H_K)$ for any non-empty index set  $K \subset \{1 \ldot d\}$. A similar result is shown for $\lambda(Q,H)$. In our notation $Q_K$ denotes the marginal df of $Q$ with respect to $K$ and $\abs{K}$ stands for the number of the elements of the index set $K$. Below $\mu_n$ and $\lambda_n$ are as defined in the Introduction.

\BT\label{GER} If $H$ is a max-id.\ df, then we have 
\bqn{\label{corWN}
	\ee = \limit{n} \mu_n(H^{1/n}, Q), \quad \lambda(Q,H)= \limit{n} \lambda_n(Q, H^{1/n}).
}
Moreover,  
\bqn{\ee &=& d+ \sum_{2 \le i\le d} (-1)^{i+1} 
	\sum_{K \subset \{1 \ldot d\},\abs{K}=i}\lambda(Q_K, H_K)
	\label{st1}
}
and 
\bqn{\lambda(Q,H) &=& d+ \sum_{2 \le i\le d} (-1)^{i+1} 
	\sum_{K \subset \{1 \ldot d\},\abs{K}=i} \mu(H_K,Q_K) \label{st2}.
}
\ET

\BRM i) 
For $H$ a max-stable df and $Q=H$ the claim in \eqref{st2} is shown in \cite{StoevDOMB}[Theorem 2.2, Eq. (13)]. \\
ii) A direct consequence of \eqref{st2} is that we can define $\lambda(Q,H)$ even if $H$ is not a 
max-id.\ df by simply using the definition of $\mu(H_K,Q_K)$. \\
iii) It is clear that $\ee \ge \mu(H_K, Q_K)$ for any non-empty index set $K\subset \{1 \ldot d\}$. Note that \eqref{corWN} shows  that exactly the opposite relation holds for $\lambda(Q,H)$ when $H$ is a max-id.\ df, namely 
$$ \lambda(Q,H) \le \lambda(Q_K, H_K).$$
\ERM 

In fact,  \eqref{st2} shows that  we can calculate both $\ee$ and $\lambda(Q,H)$ by a limit procedure if we assume that $H$ is  a max-id.\ df, see 
for more details \eqref{spiez}. Although such a limit procedure shows how to interpret these dependence  measures in terms of domination of random vectors, it does not give a precise relation with extremal properties of random samples. Therefore in the following we shall restrict our attention to the tractable case that  $H$ is a max-stable df.

\subsection{Max-stable df $H$}
We show next  the relation of $\ee$ and $\lambda(Q,H)$ with the marginal and complete domination of sample maxima mentioned in the Introduction.  Recall that in our notation $\bar F_i, \bar G_i,i\le d $ stand for the marginal survival functions of $F$ and $G$, respectively. 

\BT\label{thm20}  If $H,Q$ are max-stable df's with unit \FRE marginals  and $F,G$ are two $d$-dimensional continuous df's such that  $\limit{x} \bar F_i(x)/ \bar G_i(x)=c_i \in (0,\IF)$  for  $i\le d$ and further $F \in MDA(H), G \in MDA(Q)$, then \eqref{eqJamba} and  \eqref{lamb} hold.
\ET

\BRM The relation $\limit{n}\lambda_n(F^n, F)=\lambda(H,H)$ for $F\in MDA(H)$ is known from works of A. Gnedin, see e.g., \cite{Gnedin94,Gnedin98}. Explicit formulas are given in \cite{Ledford} for $d=2$. See also the recent contributions \cite{StoevDOMB,DombExtrC}. 
\ERM
In view of \cite{deHaan} (recall $H$ has unit \FRE marginal df's) the assumption that $H$ is max-stable implies the  following de Haan representation (see e.g., \cite{DeM15, MolchKristin}) 
\bqn{\label{deHR}
	- \ln H(x_1 \ldot x_d)= \E[\Big]{ \max_{1 \le i \le d}  \frac{Y_i}{x_i}}, \quad (x_1 \ldot x_d) \in (0,\IF)^d,
}
where $Y_j$'s are non-negative with $\E{Y_i}=1, 1 \le i \le d$.  As shown in \cite{Htilt}, see also \cite{KDEH2, Hrvoje} we have the alternative formula 
\bqn{ \label{iam}
	- \ln H(x_1 \ldot x_d)=  \sum_{i=1}^d \frac{1}{x_i} \Psi_i(x_1 \ldot x_d), \quad (x_1 \ldot x_d) \in (0,\IF)^d,
} 
where $\Psi_i$'s are non-negative  zero-homogeneous, i.e., 
$\Psi_i(cx_1 \ldot c x_d)=\Psi_i(x_1 \ldot x_d)$ for any $c>0, x_i \in (0,\IF), i\le d.$ Moreover, $\Psi_i$'s are 
bounded by 1, which  immediately  implies the validity of the lower bound in \eqref{ghr}.\\
 
In the literature $ - \ln H(1 \ldot 1)$ is also referred to as the extremal coefficient of $H$, denoted by $\theta(H)$, see e.g., \cite{StoevDOMB}.\\
Our next result gives alternative formulas for $\ee$ and shows that it is the extremal coefficient of the max-stable df $H^*$ defined by 
\bqn{\label{geret}
	- \ln H^*(x_1 \ldot x_d)= \E[\Big]{ \max_{1 \le i \le d}  \frac{ Y_i}{x_i Z_i}} 
	, \quad (x_1 \ldot x_d) \in (0,\IF)^d,
}
with $\vk Z $  being independent of $\vk Y=(Y_1 \ldot Y_d).$ Note that since 
$$ \E{Y_i}= \E{1/Z_i}=1, \quad i\le d$$ 
and $Y_i/Z_i$'s are non-negative, then $H^*$ has unit \FRE marginal df's and moreover also $\tilde H$ defined by 
\bqn{\label{geret2}
	- \ln \tilde H(x_1 \ldot x_d)= \E[\Big]{ \max_{1 \le i \le d}  \frac{ 1}{x_i Z_i}} 
	, \quad (x_1 \ldot x_d) \in (0,\IF)^d
}
 is a max-stable df with unit \FRE marginal df's. 

\BT \label{kleefeld} If  $H$ is a max-stable df with unit \FRE marginal df's and de Haan representation \eqref{geret} with $\vk Y$ being independent of  $\vk Z$ with df $Q$ which has unit \FRE marginal df's, then we have 
\bqn{\label{metrika}
	\ee&=&  \E[\Big]{ \max_{1 \le i \le d}  \frac{Y_i}{Z_i}}=  \sum_{i=1}^d \E[\Big]{\frac{1}{Z_i} \Psi_i(Z_1 \ldot Z_d)}, \quad (x_1 \ldot x_d) \in (0,\IF)^d
}
and 
\bqn{\label{maigret}
\lambda(Q,H)= \E[\Big]{ \min_{1 \le i \le d}  \frac{Y_i}{Z_i}}. 
}
 Moreover, with $H^*$ defined in \eqref{geret} 
\bqn{
	\label{9vj} \mu(H,Q) = \theta(H^*) \ge \max\Bigl(\theta(H), \theta(\tilde H) \Bigr) \ge 1
}
and 
\bqn{
	\label{maigret2} 
\lambda(Q,H) &\le & 
  \min \Biggl(  \E*{ \min_{1 \le i \le d} Y_i}, \E*{ \min_{1 \le i \le d} \frac{1}{Z_i} }  \Biggr) \le 1.
}
\ET

\BRM 
i) If  $Z_1= \cdots = Z_d= Z$ with $Z$ a unit \FRE random variable, then the zero-homogeneity of $\Psi_i$'s, 
\eqref{iam} and \eqref{metrika} imply that 
\bqn{\label{prsh} 
	\ee=  \sum_{i=1}^d \Psi_i(1 \ldot 1) \E[\Big]{ \frac{1}{Z} } = 
	\sum_{i=1}^d \Psi_i(1 \ldot 1) =  \E{\max_{1 \le i \le d} Y_i}=- \ln H(1 \ldot 1) \ge 1.
}
Further, by \eqref{maigret} we have $ \lambda(Q,H)= \E{ \min_{1 \le i \le d} Y_i}$.\\
ii) In view of \cite{StoevDOMB}[Theorem 2.2] (see also \cite{Ledford}[Eq. (6.9)]) for $H$ with de Haan representation \eqref{geret}   
$$ \lambda(H,H) = -\E[\Big]{ \frac 1 { \ln H(Y_1 \ldot Y_d)}}$$
holds, which together with  
\eqref{9vj} implies that 
$$ \mu(H_\IF,H_\IF)= \lambda(H_\IF,H_\IF)= 1$$
and thus the lower bound in \eqref{gabi} is sharp. We note in passing 
that  there are numerous papers where $\lambda_n(F^n,F)$ and $\lambda(H,H)$ appear, see e.g., \cite{de1985extremes,Ledford, Gnedin08,Gnedin12,DombExtrC} and the references therein. \\
iii) For common max-stable df's $H$ the spectral random vector $\vk Y $ that defines 
\eqref{deHR} is explicitly known. Consequently, 
for any given random vector $\vk Z$,  using the first expression in \eqref{metrika} and 
\eqref{maigret}, we can easily evaluate $\ee$ and $\lambda(Q,H)$ by Monte Carlo simulations, respectively.   
\ERM

\section{Applications}

In multivariate extreme value theory it is important to have conditions that show if a given max-stable df 
$H$ is equal to $H_0$. In case $d=2$ it is well-known that $H=H_0$ if and only if $\lambda(H,H)=0$, see \cite{StoevDOMB}[Proposition 2.2] or \cite{Gnedin94}[Theorem 2]. Consequently, when $d>2$,  
in view of   \cite{Faletal2010}[Theorem 4.3.3] we have that $H=H_0$ if and only if 
\bqn{ \lambda(H_K,H_K)=0
}
 for any index set $K\subset \{1 \ldot d\}$ with two elements. 
Therefore, in the sequel we consider for simplicity the case $d=2$ discussing some tractable conditions that are equivalent with $H=H_0$ and \eqref{turbo}. 

 As in Balkema and Resnick \cite{BaR1977}, for a given bivariate df $H$ with unit \FRE margins  define $\xi_{H} : \Rp^2 \to [0,1]$ by (set $ A=H(x_1,x_2),B= H(x_1+h, x_2+h)$)
 \BQN \label{thh}
 \quad \quad \quad \xi_{H}(\x) =\lim_{h \to 0 } \frac{[B- H((\x)+(h,-h))][B- H((\x)+(-h,h))]}{A[A+B- H((\x)+(h,-h))-H((\x)+(-h,h))]   }, \quad (\x)\in \ir.
 \EQN
  If $H$ is a continuous max-id.\ df, then in view of \cite{BaR1977} the function $\xi_H$ is non-negative, measurable and bounded by 1, almost everywhere   with respect to $dH$.\\

\BS \label{thm1} Let  $H$ and $Q$ be two biavairate df's with unit \FRE marginals. If   $H$ is a max-id.\ df, then we have 
\begin{eqnarray}\label{eq:dep:korr:1:res:1}
\lambda(Q,H) &=& \int_{\ir} [1- \xi_H(\x)] \frac{Q(\x)}{H(\x)}\, d H(\x).  
\end{eqnarray}
Moreover, if $H$ is a max-stable df, then the following conditions are equivalent: 
\begin{itemize}
	\item[i)] $H=H_0$;
	\item[ii)] $\theta(H)=- \ln H(1,1)=2$; 
		\item[iii)]	$\ee=2- \lambda(Q,H)$;   
	\item[iv)]	$\xi_H$ equals 1 almost everywhere  $dH$; 
		\item[v)] $\frac{d H^t}{ d H}= \frac{t^2 H^t}{H}$ almost everywhere  $dH$ for any $t>0$.
\end{itemize}
\label{prop1}
\ES

\BRM i) By \cite{Gnedin94}[Theorem 2] we have that 
$\lambda(H,H)=0$ is equivalent with $H=H_0$ and $\lambda(H,H)=1$ is equivalent with $H=H_\IF$. \\
ii) Statement iii) above holds for any df $Q$ with continuous marginal df's, and thus $\ee$ and $\lambda(Q,H)$ are both dependence measures for $H$.
\ERM

We conclude this section with   equivalent conditions to \eqref{turbo}. 

\BS \label{Thm2} Let $F,G$ be two continuous  bivariate df's with marginal df's $F_i,G_i,i=1,2$  satisfying $\limit{t} \bar F_i(t)/\bar G_i(t)=1$. If further $F_1,F_2$ satisfy \eqref{nco}  
 and 
$(X_1,X_2)$ 
has df $F$, then 
the following are equivalent:
\begin{itemize}
	\item[i)] $F$ has asymptotically independent components;
	\item[ii)]	$\limit{n} n\pk{X_1> n, X_2> n}=0$; 
	\item[iii)] $ \limit{n} \lambda_n(G^n,F)=0$; 
		\item[iv)]	$\limit{n}  \mu_n(F,G^n)=2$;  
	\item[v)]  $\limit{n} n \pk{ G(X_1,X_2)> 1- 1/n} =0$.
\end{itemize}
\ES 

\BRM 
i) The equivalence of i) and ii) in \neprop{Thm2} is well-known and relates to Takahashi theorem, i.e., it is enough to know that the limiting max-stable df $H$ is a product df at one point, say (1,1). See for more details in the $d$-dimensional setup \cite{Faletal2010}[p. 452].\\
ii) Recall that the assumption  $F_i \in MDA(\Phi)$ means that $\limit{n}F_i^n(a_{ni} x) = \Phi(x), x\inr$ for some norming constants $a_{ni}>0,n \in \mathbb{N}$. For notational simplicity,   in this paper  we  assume that $a_{ni}$'s equal $n$. If this is not the case, then we need to re-formulate statement ii) in \neprop{Thm2}  as
$n\limit{n}n \pk{X_1> a_{n1}, X_2 > a_{n2}}=0$. Note that if $F\in MDA(H)$ with $H$ a max-stable df, then 
\bqn{\label{coef}
	  \limit{n}n \pk{X_1> a_{n1}, X_2 > a_{n2}} = 2+ \ln H(1,1)= 2- \theta(H)=:\lambda_{F}.
}
In the literature, $\lambda_{F}$ is commonly referred to as the coefficient of upper tail dependence of $F$, see \cite{Faletal2010} for more details. 
\ERM

\section{Auxiliary Results} 
\BEL \label{sat:DEP:a}
	Let $(V_1 \ldot V_d)$ be a  random vector with continuous marginal df's $H_i,i\le d$. If further  $G$ is a $d$-dimensional df with  $G(x_1 \ldot x_d)< 1$ for any $(x_1 \ldot x_d )\in (0,\IF)^d$ and the upper endpoint of $H_i, 1\le i\le d$ equals $\IF$, then we have 
	\begin{eqnarray}\label{eq:suff:In:0}
	\limit{n} n\E{G^{n-1}(V_1\ldot  V_d)}&=& \limit{n} n \pk[\Big]{G(V_1 \ldot V_d) > 1-\frac{1}{n}} =\kappa \in [0,\IF)
	\end{eqnarray}
	if either of the limits exists. Further if 
\begin{eqnarray}\label{eq:lem:DEP:a}
G(x_1 \ldot x_d)&\le & \min_{1 \le i \le d} H_i(x_i) , \quad (x_1 \ldot x_d) \in (0,\IF)^d,
\end{eqnarray}
then  $\kappa \in [0,1].$
\EEL

\prooflem{sat:DEP:a} The proof of \eqref{eq:suff:In:0} 
follows from \cite{EHEXT}[Lemma 2.4],  see also \cite{Gnedin94}[Proposition 4]. Assuming \eqref{eq:lem:DEP:a}, if $H$ denotes the df of $(V_1 \ldot V_d)$, then  we have 
\BQNY
0 &\le & n\E{ G^{n-1}(V_1 \ldot V_d)} \le n\int_{(0,\IF)^d} \min_{1 \le i \le d} H_i^{n-1}(x_i)  \, d H(x_1 \ldot x_d)\\
&\le & n\int_{0}^\IF  H_1^{n-1}(x_1) \, d H_1(x_1)=1
\EQNY
establishing the proof. \QED 
 
 \begin{sat}\label{sat:dep:0}
	Let $F_n,G_n,n\ge 1$ be   two continuous  df's  on $[0,\IF)^d$ satisfying 
	\begin{align}\label{eq:Dini:DEP:a}
	\limit{n} F^n_n(x_1 \ldot x_d)= H(x_1 \ldot x_d), \quad \limit{n} G_n^n(x_1 \ldot x_d) = Q(x_1 \ldot x_d), \quad (x_1 \ldot x_d) \in [0,\IF)^d ,
	\end{align}
	with $H,Q$ two max-id.\ df's with unit \FRE marginal df's $\Phi$.  If for all $n$ large and some $C_1>0$ 
	\begin{eqnarray}\label{eq:dep:theo:2:con:1}
	G_n^{n}(x_1 \ldot x_d) \le C_1 \sum_{1 \le i \le d} F^{n}_{ni}(x_i) , \quad (x_1 \ldot x_d) \in (0,\IF)^d,
	\end{eqnarray}
	where $F_{ni}$ is the $i$th marginal df of $F_n$, then
	\begin{eqnarray}\label{eq:dep:korr:2:res:1}
	\limit{n} n \int_{[0,\IF)^d } G^{n}_n(x_1 \ldot x_d) \, d F_n(x_1 \ldot x_d) &=& \int_{(0,\IF)^d} Q(x_1 \ldot x_d) d\upsilon(x_1 \ldot x_d),
	\end{eqnarray}
	where $ \upsilon(\cdot)$ is the exponent measure pertaining to $H$ defined on $E:=[0, \IF]^d \setminus  \{ (0 \ldot 0) \}$.	
	If further  for all $n$ large and any $x_1
	\ldot x_d $ positive 
	\begin{eqnarray}\label{eq:dep:theo:2:con:2}
	1- G_n(x_1 \ldot x_d) \le C_2  \sum_{1 \le i \le d} \bar F_{ni}(x_i),
 	\end{eqnarray}
 	then we have 
 	%
	\begin{eqnarray}\label{eq:dep:korr:2:res:2}
	\limit{n} n \int_{[0,\IF)^d} [1- G_n(x_1 \ldot x_d )] \, d F_n^n(x_1 \ldot x_d) = -\int_{(0,\IF)^d} \ln Q(x_1 \ldot x_d) \, d H(x_1 \ldot x_d).
	\end{eqnarray}
\end{sat}

\newcommand{\ntoi}{n \to \infty }

\proofprop{sat:dep:0} 
 For notational simplicity we consider below only the case $d=2$. From the assumptions
\BQN
\label{donn}
\limit{n} F_n^n(x_{n1}, x_{n2}) =H(\x), \quad \limit{n} G_n^n(x_{n1}, x_{n2}) =Q(\x)
\EQN 
for every sequence $(x_{n1}, x_{n2}) \to (\x)\in \Rp^2$ as $n\to \IF$. 

Let $\upsilon$ be the exponent measure of $H$ defined on $E$, see \cite{Res1987} for details. For any $x_0,y_0$ positive, since by our assumptions 
$$ 
\limit{n} n[1- F_n(\x)] =- \ln H(\x)  $$
holds locally uniformly for $(\x) \in (0,\IF)^2$, using 
further  \eqref{donn} and \cite{Htilt}[Lemma 9.3] we obtain   
$$ \limit{n} 
\int_{[x_0,\IF)\times [y_0, \IF)} G^{n}_n(x_1,x_2) \, d (nF_n(x_1,x_2)) = 
\int_{[x_0,\IF)\times [y_0, \IF)} Q(x_1,x_2) \, d \upsilon(\x)=:I(x_0,y_0).
$$ 
Moreover, by \eqref{eq:dep:theo:2:con:1} 
\bqny{
\lefteqn{	n\int_{[0,\IF)^2 } G^{n}_n(x_1,x_2) \, d F_n(x_1,x_2) }\HR
	&\le & n C_1 \Bigl(  \int_{[0,x_0]} F^{n-1}_{n1}(x)\, d  F_{n1}(x)
 +  \int_{[0,y_0]} F^{n-1}_{n2}(x)\, d  F_{n2}(x)\Bigr)\HR 
	&& +\int_{[x_0,\IF)\times [y_0, \IF)} G^{n}_n(x_1,x_2) \, d (nF_n(x_1,x_2)) \\
&= & 	C_1( F_{n1}^n(x_0) + F_{n2}^n (y_0)) +\int_{[x_0,\IF)\times [y_0, \IF)} G^{n}_n(x_1,x_2) \, d (nF_n(x_1,x_2)) \\
& \to & C_1(e^{-1/x_0}+ e^{-1/y_0}) + \int_{[x_0,\IF)\times [y_0, \IF) } Q(x_1,x_2) \, d \upsilon(\x), \quad n\to \IF\HR
& \to & \int_{(0,\IF)^2 } Q(x_1,x_2) \, d \upsilon(\x), \quad x_0\downarrow 0, y_0\downarrow 0,
}
where the equality above is a consequence of the assumption that $F_n,G_n$ have continuous marginal df's. Hence \eqref{eq:dep:korr:2:res:1} follows and we show next \eqref{eq:dep:korr:2:res:2}.  Similarly, for  $x_0,y_0$ as above 
\bqny{
\lefteqn{\limsup_{\ntoi} \int_{\Rp^2} n[1-G_n(\x )] \, dF^{n}_n(\x) }\HR
	& = & \limsup_{\ntoi}\biggl[ \int_{([x_0,\IF)\times [y_0, \IF))^c}  n[1-G_n(\x )] \, dF^{n}_n(\x) \HR 
	&&  +\int_{[x_0,\IF)\times [y_0, \IF)} n[1-G_n(\x )] \,dF^{n}_n(\x)\biggr]
	\HR
	& \le & C_2\limsup_{\ntoi} \int_{([x_0,\IF)\times [y_0, \IF))^c} n[\bar F_{n1}(x_1)+ \bar F_{n2}(x_2)]  \,dF^{n}_n(\x)
	\HR
	&& +\limsup_{\ntoi}\int_{[x_0,\IF)\times [y_0, \IF)} n[1-G_n(\x )] \,dF^{n}(\x)
	\HR
	& \le  &C_2\limsup_{\ntoi}(F_{n1}^{n}( x_0 )+ F_{n2}^{n}( y_0) )  \Bigl[ n \bar F_{n1}(x_0) +n \bar F_{n2}( y_0) \Bigr] \HR
&& - \int_{[x_0,\IF)\times [y_0, \IF)}\ln Q(\x) \, d H(\x)\HR
& =&  C_2\Bigl[ e^{-1/x_0}+ e^{- 1/y_0} \Bigr]\Bigl [ \frac{1}{x_0}+ \frac{1}{y_0} \Bigr]
- \int_{[x_0,\IF)\times [y_0, \IF)}\ln Q(\x) \, d H(\x) \HR
&\to & -\int_{(0,\IF)^2 }\ln Q(\x) \, d H(\x) , \quad x_0 \downarrow 0, y_0 \downarrow  0,
} 
hence the proof follows. \QED

\BRM The validity of \eqref{eq:dep:theo:2:con:1} has been shown under the  assumption that $G_n$ is a continuous df. From the proof above it is easy to see that \eqref{eq:dep:theo:2:con:1} still holds if we assume instead that  $G_n$ is continuous and positive such that $G_n^n$ is a df. Similarly, for the validity of \eqref{eq:dep:korr:2:res:2} it is enough to assume that $F_n^{n}$ is a continuous  df.  
\label{rma}
\ERM

{
\begin{korr} \label{kLeje} If $H$ is  a bivariate max-stable df with unit \FRE marginal df's $H_1$ and $H_2$, then for $u,t$ positive
	\BQN
	\int_{\ir} \min \bigl(H_1^{1/u}(x_1),\, H^{1/t}_2(x_2) \bigr)\dln&=&
	u+t+ \ln H(1/u,1/t).
	\EQN
\end{korr}
\proofkorr{kLeje} The proof follows using Fubini Theorem and the homogeneity property of the exponent measure inherited by \eqref{hany}. We give below an alternative proof. Let  $(V_1,V_2)$ have df $H$ and set $U_i=H_i(V_i),i=1,2$. By the assumptions since the df $H$ is continuous, applying \netheo{thm20} and \eqref{eq:suff:In:0}  with $u,t>0$ we obtain 
\begin{eqnarray*}
	\lefteqn{\int_{\ir} \min \bigl( H^{1/u}_1(x_1 ), H_2^{1/t}( x_2) \bigr) \dln}
	\HR
	&=& \limit{n}    n \int_{(0,\IF)^2} \min \bigl( H^{n/u}_1(x_1),  H_2^{n/t}(x_2) \bigr) \, d H(\x)
	\HR
	&=& \limit{n}    n \pk[\Big]{ \min \bigl( H^{1/u}_1(V_1), H_2^{1/t}(V_2) \bigr)>1-\frac{1}{n}}
	\HR
	&=& \limit{n} n \pk[\Big] { U_1>1-\frac{u}{n}, U_2 >1-\frac{t}{n}}
=	u+t+ \ln H(1/u,1/t)
\end{eqnarray*}
establishing the proof. \QED 
\COM{Alternatively, since $ H^s(s \x)=H(\x)$ for any $\x \in (0,\IF)^2$ and thus 
$ \upsilon\{sB\}=s^{-1} \upsilon\{B\}, $
with $B$ Borel sets bounded away from origin, by Fubini-Tonelli theorem (take for simplicity $u=t=1)$)
\BQNY
\int_{\Rp^2} \min( e^{- 1/x_1}, e^{-1/x_2}) \, \upsilon(d \x) &=& \int_{0}^1 \upsilon \{(\x)\in (0,\IF)^2:  \min( e^{- 1/x_1}, e^{-1/x_2} )>s \} \, ds
\HR
&=&  \int_{0}^1 \upsilon \{ ( -1/ \ln s , \infty ]^2\} \, ds
\HR
&=& \upsilon\{(1, \IF)^2\} \int_{0}^1 ( -\ln s) \, ds
\HR
&=&\limit{n}n  \pk[\Big]{H_1(V_1)>1-\frac{1}{n},H_2(V_2)>1-\frac{1}{n}}
\HR
&=&2+ \ln H(1,1),
\EQNY
where the second last step follows by \nekorr{K1} taking  
$$F(\x)=H(\x), \quad G(\x)=H_1(x_1)H_2(x_2).$$
   \QED
}

\section{Proofs}
\prooftheo{GER} For $n>0$ set   $A_n=Q^{1/n}$ and $B_n= H^{1/n}$. Since $H$ is a max-id.\ df, then $B_n$ is a df for any $n>0$. Furthermore, since $H_i=Q_i, i\le d$ (recall $H_i,Q_i$ are the marginal df's of $H$ and $Q$, respectively),  it can be easily checked that 
we can apply \neprop{sat:dep:0}, which together with  Remark \ref{rma}  imply 
\bqn{ \label{spiez} \lefteqn{\limit{n} n \int_{\R^d  } [1- H^{1/n} (x_1 \ldot x_d)]\, d Q(x_1 \ldot x_d)}\notag \\
	&=&  \limit{n} n \int_{\R^d  } [1- B_n  (x_1 \ldot x_d)]\, d A_n^n(x_1 \ldot x_d)\notag \\
&=& - \int_{\R^d}  \ln H (x_1 \ldot x_d)]\, d Q(x_1 \ldot x_d) =\mu(H,Q).
}
The second claim in \eqref{corWN} follows with similar arguments and  therefore we omit its proof. \\
Next, for any non-empty subset $K$ of $\{1 \ldot d\}$  with  $m=\abs{K}$ elements  by  \eqref{corWN}  
$$ \mu(H_K,Q_K) = \limit{n} n \int_{\R^{m }} [1- F_{n,K}(x_1 \ldot x_d)]\, d Q_K(x_1 \ldot x_d)$$
and 
$$ \lambda(Q_K,H_K)=  \limit{n} n \int_{\R^{m}} Q_K (x_1 \ldot x_d)\, d F_{nK}(x_1 \ldot x_d),$$
where  $F_{nK}, Q_K$ are the marginals of $F_n$ and $Q$ with respect to $K$. Note that for notational simplicity we write the marginal df's with respect to $K$ as functions of $x_1 \ldot x_d$ and not as functions of $x_{j1} \ldot x_{jm}$ where $K= \{j_1 \ldot j_m\}$ has $m=\abs{K}$ elements.
By Fubini Theorem 
$$ \int_{\R^{m}} Q_K (x_1 \ldot x_d)\, d F_{nK}(x_1 \ldot x_d) = 
\int_{\R^{m}} \overline{F}_{nK} (x_1 \ldot x_d)\, d Q_K(x_1 \ldot x_d),$$
where $\overline{F}_{nK}$ stands for the joint survival function of $F_{nK}$. In the light of the  inclusion-exclusion formula 
$$ 1- F_{n}(x_1 \ldot  x_d) = \sum_{1 \le i\le d} (-1)^{i+1} 
\sum_{K \subset \{1 \ldot d\},\abs{K}=i}\overline{F}_{nK}(x_1 \ldot x_d), \quad (x_1 \ldot x_d) \in \R^d.
$$
Using further the fact that $H$ and $Q$ have the same marginal df's, for any index set $K$ with only one element we have 
$$ \limit{n} n\int_{\R^d} \overline{H}_{nK}(x_1 \ldot x_d)\, d Q(x_1 \ldot x_d) = 
\limit{n} n \int_0^1 (1- t^{1/n}) d t = 1,$$
hence 
\bqny{ \lefteqn{\ee}\\ 
	&=& 
	\limit{n} n \int_{\R^{\abs{K}}} [1- F_{n}(x_1 \ldot x_d)]\, d Q(x_1 \ldot x_d)\\
&=& 	d+ \limit{n} n\int_{\R^{\abs{K}}} \sum_{2 \le i\le d} (-1)^{i+1} 
\sum_{K \subset \{1 \ldot d\},\abs{K}=i}\overline{F}_{nK}(x_1 \ldot x_d)\, d Q(x_1 \ldot x_d)\\
&=&  d+  \sum_{2 \le i\le d} (-1)^{i+1}\limit{n}   n\int_{\R^{\abs{K}}}
\sum_{K \subset \{1 \ldot d\},\abs{K}=i}\overline{F}_{nK}(x_1 \ldot x_d) \, d Q_K (x_1 \ldot x_d)\\
&=&  d + \sum_{2 \le i\le d} (-1)^{i+1} 
\sum_{K \subset \{1 \ldot d\},\abs{K}=i} \lambda(Q_K, H_K)
}
and thus \eqref{st1} follows. Since by the inclusion-exclusion formula we have further 
$$ \overline{F}_n(x_1 \ldot  x_d) = \sum_{1 \le i\le d} (-1)^{i+1} 
\sum_{K \subset \{1 \ldot d\},\abs{K}=i}[1-{F}_{nK}(x_1 \ldot x_d)], \quad (x_1 \ldot x_d) \in \R^d
$$
the claim in \eqref{st2} follows with similar arguments as above.  
\QED

\prooftheo{kleefeld} The claim in \eqref{metrika} follows by the de Haan and inf-argmax representation of $H$. Since by the independence of $Y_i$'s and $Z_i$'s and the fact that 
$\E{Y_i}=\E{1/Z_i}=1$  we have that 
\bqn{
	\E{ Y_i/Z_i}=\E{Y_i}\E{1/Z_i} =1
\label{liebes}
}
is valid for any $i\le d$. Consequently, by \eqref{st2}, \eqref{metrika} and the fact that for given constants 
$c_1 \ldot c_d$ 
$$ \min _{1 \le i \le d} c_i =   \sum_{i=1}^d (-1)^{i+1}  \sum_{ K \subset  \{ 1 \ldot d \}: \abs{K}=i} \max_{j \in K } c_j,
$$
then we have 
\bqny{\lambda(Q,H) &=& \sum_{i=1}^d \E*{\frac{Y_i}{Z_i}} + 
	\sum_{2 \le i\le d} (-1)^{i+1} 
	\sum_{K \subset \{1 \ldot d\},\abs{K}=i} \E*{ \max_{ j \in K}  \frac{ {  Y_j}}{Z_j}}\\
	&=& \E*{ \min_{ 1 \le j \le d }  \frac{ {  Y_j}}{Z_j}}
}
establishing \eqref{maigret}.\\
Further, since  \eqref{liebes} holds, then  by de Haan representation of max-stable df's  we have that the df's $H^*, \tilde H$ defined  in \eqref{geret} and \eqref{geret2}, respectively 
are max-stable with unit \FRE marginal df's. Hence \eqref{metrika} implies that 
$\ee= \theta(H^*)$. Note in passing that for $Q=H$ this follows also from \cite{StoevDOMB}[Proposition 2.2]. \\
Using again that $Y_i$'s are independent of $Z_i$'s and $\E{Y_i}=1, i\le d$ we obtain  (recall $Y_i$'s and $Z_i$'s are non-negative random variables)
\bqny{
	\ee &=& \E[\Big]{ \E[\Big]{ \max_{1 \le i \le d}  \frac{ Y_i}{ Z_i} \Bigl \lvert (Z_1 \ldot Z_d)}}\\
	&\ge & \E[\Big]{ \max_{1 \le i \le d}  \frac{ \E{  Y_i}}{Z_i}}\\ 
	& \ge & \E[\Big]{ \max_{1\le i \le d} \frac{1}{Z_i}} = \theta(\tilde H)\\ 
	& \ge & \max_{1 \le i \le d } \E[\Big]{ \frac{1}{Z_i}}=1.
}
With the same arguments  using now that $\E{1/Z_i}=1, i\le d$ we have 
\bqny{
	\ee   &=& \E[\Big]{ \E[\Big]{ \max_{1 \le i \le d}  \frac{ Y_i}{ Z_i} \Bigl \lvert (Y_1 \ldot Y_d)}}\\
	&\ge &  \E[\Big]{ \max_{1\le i \le d} Y_i }  = - \ln H(1 \ldot 1) = \theta(H).
}
The lower bound in \eqref{maigret2} follows with similar arguments, hence the proof is complete.  
\COM{\bqny{
	\lambda(Q,H)   &=& \E[\Big]{ \E[\Big]{ \min_{1 \le i \le d}  \frac{ Y_i}{ Z_i} \Bigl \lvert (Y_1 \ldot Y_d)}}\\
	&\le &  \E[\Big]{ \min_{1\le i \le d} Y_i } 
}
\bqny{
	\lambda(Q,H)   &=& \E[\Big]{ \E[\Big]{ \min_{1 \le i \le d}  \frac{ Y_i}{ Z_i} \Bigl \lvert (Y_1 \ldot Y_d)}}\\
	&\le &  \E[\Big]{ \min_{1\le i \le d}  \frac 1 {Z_i} } 
}
}
\QED

\prooftheo{thm20} Suppose without loss of generality that $F$ satisfies \eqref{plakut}. If $F_i=G_i,i=1,2$, then the claim follows from  \nelem{sat:DEP:a} and \neprop{sat:dep:0}.  We consider next the general case that $F_i$'s are tail equivalent to $G_i$'s and suppose for simplicity that $d=2$. In view of \cite{EHEXT}[Lemma 2.4] we have 
$$ \limit{n} n \int_{[0,\IF)} G_i^n(x) d F_i(x) = c_i \in [0,\IF), \quad i=1,2$$
if and only if $\limit{n} n \pk{ G_i(X_i)> 1-1/n}= c_i$ or equivalently 
$$ \limit{x} \frac{\bar F_i(x)}{\bar G_i(x)}=c_i.$$
By the assumption $c_i \in (0,\IF)$ for $i=1,2$.
Consequently, for all $x>0$ there exist $a_1,a_2$ positive such that 
$$a_1 \bar F_i(x) \le \bar G_i(x) \le  a_2 \bar F_i(x).$$
Assume for simplicity that $c_i=1,i=1,2$. By the assumptions   
$$ n \bar F_i (nx) \to 1/x, \quad n \bar G_i (nx) \to 1/x , \quad n\to \IF$$
uniformly for $x$ in  $[t,\IF), t>0$. Further,  for $i=1,2$ we have 
$$ \lim_{t\downarrow  0} \limit{n} n\int_{[0,t]} G_i^n(nx) \, d F_i(nx)= \lim_{t\downarrow  0} \limit{n} n\int_{[0,t]} \bar G_i(nx) \, d F_i^n(nx)=0,$$
which implies 
$$ \lim_{t\downarrow  0} \limit{n} n\int_{[0,t]^2} G^n(nx,ny) \, d F(nx,ny)
= \lim_{t\downarrow  0} \limit{n} n\int_{[0,t]^2} [1- G(nx,ny)] \, d F^n(nx,ny)=0.$$
As in the proof of \neprop{sat:dep:0},  using that $F$ and $G$ are in the MDA of $H$ and $Q$, respectively, it follows that for any integer $k$
$$ \limit{n} n\int_{[0,\IF)^2} G^{n-k}(\x) \, d F(\x)= 
\int_{(0,\IF)^2} Q(\x) \, d \upsilon(\x)= \lambda(Q,H)$$
and further 
$$ \limit{n} n\int_{[0,\IF)^2} [1- F(\x)] \, d G^{n-k}(\x)= 
- \int_{(0,\IF)^2} \ln H(\x) \, d Q(\x)= \ee $$
establishing the proof.  
\QED

\proofprop{thm1}  In view of \netheo{GER}, since $H$ being a max-id.\ df implies that $H^{1/n}$ is a  df for any $n\ge 1$  we have with $F_n= Q^{1/n}$ 
\bqny{
	\int_{(0,\IF)^2} Q(\x) d \upsilon(\x)&=&\limit{n} n\int_{\ir} Q(\x) \, d H^{1/n} (\x) \HR
	&=& 2-\limit{n} n \int_{[0,\IF)^d } [1- H^{1/n}(\x)] \, d F_n^n (\x) \\
	&=& \int_{\ir} [2 + \ln H(\x)] \, d Q(\x).  
}
Since further by \cite{BaR1977}[Theorem 7] the restriction of $\upsilon$ on $(0,\IF)^2$ denoted by $\upsilon_0$ satisfies  
$$ \dfrac{d \upsilon_0  }{d H } =\frac{1-\xi_H}{H}$$
and $\xi_H(\x) \in [0,1]$ almost everywhere  $dH$, 
then the first claim  follows.

 The equivalence of i) and ii) is known as Takahashi Theorem, see \cite{Faletal2010}[Theorem 4.3.2]. Since $\xi_H \in [0,1]$ almost everywhere $dH$, then  the equivalence of ii) and iii) is a direct consequence of \eqref{eq:dep:korr:1:res:1} and the fact that 
 $\lambda(Q,H)= 2- \mu(H,Q)$, see \eqref{st2}. Clearly, by \eqref{eq:dep:korr:1:res:1} we have thus $\xi_H=1$ almost everywhere   $dH$ is equivalent with $H=H_0$, whereas iv) is equivalent with v) is consequence of \cite{BaR1977}[Theorem 7].  \QED 
 
\proofprop{Thm2} 
 If $F$ \eqref{turbo} holds, then clearly ii) is satisfied and thus i) implies ii).  If ii) holds, then 
 \bqny{
 	 \limsup_{n\to \IF} F^n(nx_1,nx_2) &=& \exp\Bigl( 	 \limsup_{n\to \IF} n \ln (1- [1- F(nx_1,nx_2)] ) \Bigr )\\
 	 &= & \exp\Bigl(- 	 \limsup_{n\to \IF} n ( 1- F(nx_1,nx_2)) \Bigr)\\
 	 &\le  & \exp\Bigl( - 	 \limsup_{n\to \IF} [  n \bar F_1(nx_1) + n \bar F_2(nx_2)- 
 	 	n \pk{X_1> n \min(x_1,x_2), X_2> n \min(x_1,x_2) }]\Bigr)\\
 	 &=  & \exp\bigl( -1/x_1- 1/x_2 \bigr), \quad x_1,x_2>0.
 }
As for the derivation of \eqref{ghr} we obtain  further 
\bqn{ \label{bacb}
		  	 \liminf_{n\to \IF}  F^n(nx_1,nx_2) \ge  \exp\bigl( -1/x_1- 1/x_2 \bigr), \quad
 x_1,x_2>0
}
implying that $F\in MDA(H_0)$, hence i) follows. \\
Assuming iii) and since the marginal df's of $G$ are in the MDA of $\Phi$, with the same calculations as in \eqref{bacb} for the df $G$ we obtain
\bqny{ 0 = \limit{n} n \int_{[0,\IF)^2} G^n(\x) \, d F(\x) &\ge& 
	\limit{n} n \pk{X_1> n, X_2> n} G^n(n,n)\\
& \ge & c 		\limit{n} n \pk{X_1> n, X_2> n} 
}
for some $c\in (0, e^{-2})$, hence ii) follows.\\ 
Next, assume that  ii) holds.  We have that $G(\x) \le G_1(x_1)G_2(x_2)=:K(\x)$  and by the assumption that $G_i$'s are in the MDA of $\Phi$ it follows that  $K$ is in the MDA of $H_\IF$. Further ii) implies that $F\in MDA(H_0)$ and $ - \ln H(1,1)=2$. Consequently,  
\netheo{thm20} yields  
$$ \limit{n} \lambda_n(K^n, F)= \lambda(H_\IF, H_0).$$ 
But from \nekorr{kLeje} we have  that $\lambda(H_\IF, H_0)=0$, hence ii) implies iii).\\
Let $\overline G$ be the joint survival function of the bivariate df $G$. For any positive integer $n$, we have that $F^n$ is a bivariate df. Hence by Fubini theorem and the fact that $F_i=G_i,i=1,2$ are continuous df's, for any positive integer $n$ we obtain 
\bqny{  \int_{ \R^2 } F^{n}(x_1,x_d) d G(x_1,x_2)
	 &=&  
	\int_{\R^2 } \overline{G}(x_1,x_2)  d F^{n}(x_1,x_2) \\
	&=& 2 n  / (n+1) 
	- 	\int_{ \R^2 } [1- G(x_1,x_2)]  d F^{n}(x_1,x_2) 
}
and thus the equivalence of iii) and iv) follows. The equivalence  iv) and v) follows from \nelem{sat:DEP:a} and thus the proof is complete.
\QED 

	\section*{Acknowledgments}
	 I am in debt  to both reviewers for numerous suggestions and corrections. 	 Support from  SNSF Grant 200021-175752/1 is kindly acknowledged.
	\bibliographystyle{ieeetr}
	\bibliography{EEEA}
\end{document}